\documentclass[a4paper,11pt]{article}
\usepackage{amsfonts}
\usepackage{bbding}
\usepackage{amssymb}
\usepackage{mathrsfs}
\usepackage{float}
\usepackage{amsmath}
\usepackage{amsfonts,amssymb}
\usepackage{graphicx,color}
\usepackage{mathrsfs}
\usepackage{graphicx}%
\catcode`\@=11
\@addtoreset{equation}{section}

\catcode`\@=12

\newtheorem{Theorem}{Theorem}[section]

\newtheorem{Lemma}{Lemma}[section]
\newtheorem{Example}{Example}[section]

\newtheorem{Remark}{Remark}[section]

\newtheorem{Proposition}{Proposition}[section]
\newtheorem{Conjecture}{Conjecture}[section]

\def\2{{I \hskip -1.0mm I}}
\def\3{{I \hskip -1.0mm I\hskip -1.0mm I}}
\def\4{{I \hskip -0.9mm V}}
\def\6{{V \hskip -1.35mm I}}

\begin{document}
\title{A new kind of functional differential equations}
\author{De-Xing Kong\footnote{Corresponding author: dkong@zju.edu.cn} $\quad$ and $\quad$ Cheng Zhang\\ \\
{\it {\small Department of Mathematics, Zhejiang University}}\\
{\it  {\small Hangzhou 310027, China}}}
\date{}
\maketitle

\begin{abstract}
     In this paper we introduce and investigate a new kind of functional (including ordinary and evolutionary partial) differential equations. The main goal of this paper is to explore our new philosophy by some examples on functional ODEs and PDEs. For some typical examples, we prove the global existence of smooth solutions, analyze some interesting properties enjoyed by these solutions, and illustrate the differences between this new class of equations and the traditional ones. This kind of functional differential equations is a new and powerful tool to study some problems arising from both mathematics and physics, more applications in particular to differential geometry and fundamental physics can be expected.

\vskip 6mm


\vskip 3mm


\end{abstract}

\newpage\baselineskip=7mm

\section{Introduction}

A differential equation is an equation involving an unknown function of one or two or more variables and certain its (partial) derivatives. In general, we can write out symbolically a typical differential equation, as follows. Fix an integer $k\ge 1$ and let $\Omega$ denote an open subset of $\mathbb{R}^n\;(n\ge 0)$. An expression of the form
\begin{equation}\label{1.1}
\mathscr{P}(D^ku(t,x),D^{k-1}u(t,x),\cdots,Du(t,x),u(t,x),t,x)=0\quad (x=(x_1\cdots,x_n)\in \Omega)
\end{equation}
is called a $k$-th order differential equation, where
$$\mathscr{P}:\;\;\mathbb{R}^{(n+1)^k}\times\mathbb{R}^{(n+1)^{k-1}}\times\cdots\times\mathbb{R}^{n+1}\times\mathbb{R}\times\mathbb{R}\times\Omega\rightarrow \mathbb{R}$$
is given and $u=u(t,x_1\cdots,x_n)$ is the unknown. The variable $t$ always denotes time. Obviously, if $n=0$, then (\ref{1.1}) is an ordinary differential equation; if $n\ge 1$, then it is partial differential equations.

Instead of (\ref{1.1}), we are interested in the following functional differential equation
\begin{equation}\label{1.2}
\mathscr{F}(D^ku(t,x),D^ku(t/2,x),D^{k-1}u(t,x),D^{k-1}u(t/2,x),\cdots,Du(t,x),Du(t/2,x),u(t,x),u(t/2,x),t,x)=0,
\end{equation}
where
$$\mathscr{F}:\;\;\mathbb{R}^{(n+1)^k}\times\mathbb{R}^{(n+1)^k}\times\mathbb{R}^{(n+1)^{k-1}}\times\mathbb{R}^{(n+1)^{k-1}}
\times\cdots\times\mathbb{R}^{n+1}\times\mathbb{R}^{n+1}\times\mathbb{R}\times\mathbb{R}\times\mathbb{R}\times\Omega\rightarrow \mathbb{R}$$
is given and $u=u(t,x_1\cdots,x_n)$ is the unknown. The equation (\ref{1.2}) is also known as ``nonlocal", some time called {\it time-delay}\footnote{Usually, the traditional time-delay equation means that the time variable in the equation is $t-r$ for some constant $r$ instead of $t/2$ in this paper (see \cite{m}).}, differential
equations. This kind of new time delay effect can also be introduced in the study of some problems arising in fundamental physics in a way similar to the work \cite{am}. In this paper, we shall explore our new philosophy by some simple but interesting examples.

In fact, $t/2$ in (\ref{1.2}) can be replaced by $t/\alpha \; (\alpha > 1)$ or $\phi(t)$ with $0<\phi(t)<t$. In this case, we can develop a similar theory.

Some typical examples are as follows:

\subsection{Functional ODEs}

$\;\;\;\;\;\;$ {\bf {- Linear-like equation}}
\begin{equation}\label{1.3}
\frac {dx\left(t\right)}{dt}=x\left(\frac t2\right).
\end{equation}

\vskip 3mm
{\bf {- Riccati-like equation }}
\begin{equation}\label{1.4}
\frac {dx\left(t\right)}{dt}=x^2\left(\frac t2\right).
\end{equation}

\vskip 3mm
{\bf {- A system of functional ODEs}}
\begin{equation}\label{1.5}
\frac {dx\left(t\right)}{dt}=A\left(t\right)x\left(t\right)+F\left
(x\left(\phi\left(t\right)\right)\right), \quad t\geq 0,
\end{equation}
where $x\left(t\right)=\left(x_1\left(t\right),\cdots
,x_n\left(t\right )\right)^T$, $A\left(t\right)$ is a continuous
$n\times n$ matrix-valued function of $t\in\mathbb{R}^{+}$, $F=\left(F_1,
\cdots ,F_n\right)^T$ is a locally Lipschitzian continuous
vector-valued function of $x\in\mathbb{R}^n$, and $\phi$ is a continuous
real-valued function of $t\in\mathbb{R}^{+}$ with
\begin{equation}\label{1.6}
0<\phi\left(t\right)<t\quad\text{ for all }\,\,t> 0.\end{equation}

\subsection{Functional PDEs}

$\;\;\;\;\;\;${\bf {- Burgers-like equation}}
\begin{equation}\label{1.7}
u_t(t,x)+u(t/2,x)u_x(t,x)=0.
\end{equation}

{\bf {- Heat-like equation}}
\begin{equation}\label{1.8}
u_t(t,x)- u_{xx}(t/2,x)=0.
\end{equation}

{\bf {-  Wave-like equation}}
\begin{equation}\label{1.9}
u_{tt}(t,x)- u_{xx}(t/2,x)=0.
\end{equation}

\subsection{Functional geometric flows}

$\;\;\;\;\;\;$ {\bf{- Functional Ricci flow}}

Let $\mathscr{M}$ be an $n$-dimensional complete Riemannian manifold
with Riemannian metric $g_{ij}$, the Hamilton's Ricci flow is given
by the evolution equation (see \cite{h})
\begin{equation}\label{1.10}
\frac{\partial g_{ij}}{\partial t}=-2R_{ij}
\end{equation}
for a family of Riemannian metrics $g_{ij}(t)$ on $\mathscr{M}$.
(\ref{1.10}) is a nonlinear system of second order partial differential
equations on the metric $g_{ij}$. The functional Ricci flow under
consideration reads
\begin{equation}\label{1.11}
\frac{\partial g_{ij}}{\partial t}=-2
\left.R_{ij}\right|_{\frac{t}{2}},
\end{equation}
where $\bullet|_{\frac t2}$  stands for the value of $\bullet$ at
the time $\frac t2$. Thus, in (\ref{1.11}),
$\left.R_{ij}\right|_{\frac{t}{2}}$ is the Ricci curvature
corresponding to the metric $g_{ij}(t/2)$, i.e., the metric at the
time $\frac t2$. Obviously, (\ref{1.11}) is a system of functional partial
differential equations. Noting the special time-delay effect, we may
consider the Cauchy problem for (\ref{1.11}) with the initial data
\begin{equation}
t=0:\;\;g_{ij}=g_{ij}^0,
\end{equation}\label{1.12}
where $g_{ij}^0$ is a given Riemannian metric on $\mathscr{M}$.
Therefore, the Cauchy problem (\ref{1.11})-(1.12) gives the evolution of the
metric $g_{ij}^0$ under the flow (\ref{1.11}).

\vskip 3mm

{\bf{- Functional mean curvature flow}}

Let $\mathscr{M}$ be an $n$-dimensional smooth manifold and
$$X(\cdot,t): \quad \mathscr{M} \rightarrow \mathbb{R}^{n+1}$$ be a one-parameter
family of smooth hypersuface immersions in $\mathbb{R}^{n+1}.$ The
traditional mean curvature flow is defined by (see \cite{eh}-\cite{g} and \cite{hu1})
\begin{equation}\label{1.13}
\dfrac{\partial }{\partial
t}X(x,t)=H(x,t) \vec{n}(x,t),\qquad \forall ~x\in\mathscr{M} ,\quad
 \forall ~t\geq0,\end{equation}
 where $H(x,t)$ is the mean curvature of $X(x,t)$ and $\vec n(x,t)$ is
the unit inner normal vector on $X(\cdot, t)$. The nonlocal mean
curvature flow under consideration in the present paper is
\begin{equation}\label{1.14}
\dfrac{\partial }{\partial
t}X(x,t)=H(x,t/2) \vec{n}(x,t/2),\qquad \forall ~x\in\mathscr{M}
,\quad
 \forall ~t\geq0,\end{equation}
It is easy to see that (\ref{1.14}) is a system of nonlocal partial
differential equations. Noting the special time-delay effect, we may
consider the Cauchy problem for (\ref{1.14}) with the initial data
\begin{equation}\label{1.15}
t=0:\;\;X=X_0,
\end{equation}
where $X_0$ is a given hypersurface. Thus, the Cauchy problem
(\ref{1.15}) gives the evolution of the hypersurface $X_0$ under the
flow (\ref{1.14}).
\begin{Remark} {\it For the inverse mean curvature flow (see \cite{hi}) and the hyperbolic mean curvature flow (see \cite{hkl}, \cite{klw} and \cite{le}), we have a similar discussion.}
\end{Remark}

\vskip 3mm

{\bf{- Functional hyperbolic geometric flow}}

Let $\mathscr{M}$ be an $n$-dimensional complete Riemannian manifold
with Riemannian metric $g_{ij}$, the hyperbolic geometric flow
considered in \cite{k}-\cite{kl} and \cite{klx} is described by the evolution equation
\begin{equation}\label{1.16}
\frac{\partial^{2}g_{ij}}{\partial t^{2}}=-2R_{ij}
\end{equation}
for a family of Riemannian metrics $g_{ij}(t)$ on $\mathscr{M}$.
(\ref{1.16})  is a nonlinear system of second order partial differential
equations on the metric $g_{ij}$. The functional hyperbolic geometric
flow under considered here reads
\begin{equation}\label{1.17}
\frac{\partial^{2}g_{ij}}{\partial t^{2}}=-2
\left.R_{ij}\right|_{\frac{t}{2}}.
\end{equation}
Obviously, (\ref{1.17}) is a system of functional partial differential
equations. Noting the special time-delay effect, we may consider the
Cauchy problem for (\ref{1.17})  with the initial data
\begin{equation}\label{1.18}
t=0:\;\;g_{ij}=g_{ij}^0,\quad \frac{\partial g_{ij}}{\partial
t}=k_{ij}^0,
\end{equation}
where $g_{ij}^0$ is a given Riemannian metric on $\mathscr{M}$, and
$k^0_{ij}$ is a symmetric tensor on $\mathscr{M}$. Therefore, the
Cauchy problem (\ref{1.17})-(\ref{1.18}) gives the evolution of the metric
$g_{ij}^0$ under the flow (\ref{1.17}).

\subsection{A shifted view of fundamental physics}

Atiyah and Moore speculated on the role of
relativistic versions of delayed differential equations in
fundamental physics (see \cite{am}). Since relativistic invariance
implies that one must consider both advanced and retarded terms in
the equations, in \cite{am} they refereed to them as shifted
equations and showed that the shifted Dirac equation has some novel
properties and a tentative formulation of shifted Einstein-Maxwell
equations naturally incorporates a small but nonzero cosmological
constant.

\subsection{Aim and outline of the paper}

The main goal of this paper is to explore our new philosophy by some examples on functional ODEs and PDEs mentioned above. For typical examples\footnote{In fact, for the hyperbolic equations (\ref{1.7}) and (\ref{1.9}), we can prove the global existence of smooth solutions in a manner similar to \cite{k-li}.} (e.g., (\ref{1.3}) and (\ref{1.8})), we prove the global existence of smooth solutions, analyze some interesting properties enjoyed by these solutions, and illustrate the differences between this new class of equations and the traditional ones.

The paper is organized as follows. In Section 2, we explore our philosophy on this new kind of functional ODEs and PDEs under consideration. Sections 3 and 4 are devoted to investigating the equation (\ref{1.3}) and presenting some interesting properties enjoyed by the solution of (\ref{1.3}),  while Section 5 is devoted to studying the equation (\ref{1.8}).

\section{New philosophy}

In this section, we explore our new philosophy and the differences between functional differential equations and traditional ones by some examples.

\begin{Example}
Taking the equation (\ref{1.4}) as the first example, we consider the Cauchy problem for (\ref{1.4}), i.e.,
$$ \frac {dx\left(t\right)}{dt}=x^2\left(\frac t2\right)$$
with the initial data
\begin{equation}\label{2.1}
x(0)=1.
\end{equation}
It is easy to see that the Cauchy problem (\ref{1.4}), (\ref{2.1}) has the following global smooth solution
\begin{equation}\label{2.2}
x\left(t\right)=e^t\quad\text{  for all }\,\,t.
\end{equation}
On the other hand, the solution of the Cacuhy problem for the traditional Riccati equation
\begin{equation}\label{2.3}
 \frac {dx\left(t\right)}{dt}=x^2\left(t\right)
\end{equation}
with the initial data (\ref{2.1}) reads
\begin{equation}\label{2.4}
x\left(t\right)=\frac 1{1-t}.\end{equation}
Clearly, the solution blows up at the time $t=1$, since
\begin{equation}\label{2.5}
x\left(t\right)\nearrow\infty\quad\text{  as }\qquad t\nearrow 1
.\end{equation}
\end{Example}

\begin{Example} We next consider another simple but interesting example --- the equation (\ref{1.3}). To illustrate the difference between the equation (\ref{1.3}) and the traditional linear ODE, we first recall the solution of the following Cauchy problem
\begin{equation}\label{2.9}
\left\{\aligned
&\frac {dx\left(t\right)}{dt}=x\left(t\right),\\
&x\left(0\right)=1.\endaligned \right.\end{equation}
It is well-known that the solution can be explicitly expressed by
\begin{equation}\label{2.10}
x\left(t\right)=e^{t}\quad\text{  for all }\,\,t.\end{equation}
Obviously, the solution is a strictly increasing function for $t$.
However, for the following Cauchy problem for the equation (\ref{1.3})
\begin{equation}\label{2.11}
\left\{\aligned
&\frac {dx\left(t\right)}{dt}=x\left(\frac t2\right)\\
&x\left(0\right)=1,\endaligned ,\right.\end{equation}
in next section, we will prove

\begin{Theorem} The Cauchy problem (\ref{2.11}) admits the
global solution shown as in the following figures:

\begin{figure}[H]
    \begin{center}
\vspace{0.3cm}
\scalebox{0.8}[0.55]{\includegraphics[trim=0mm 0mm
0mm 0mm, clip, width=9cm, height=12cm]{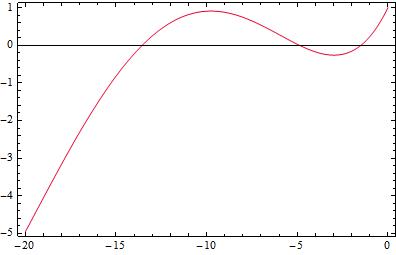}}
\caption{The figure of the solution on the interval
[-20,0]}
    \end{center}
\end{figure}

\begin{figure}[H]
    \begin{center}
\vspace{0.3cm}
\scalebox{0.8}[0.55]{\includegraphics[trim=0mm 0mm
0mm 0mm, clip, width=9cm, height=12cm]{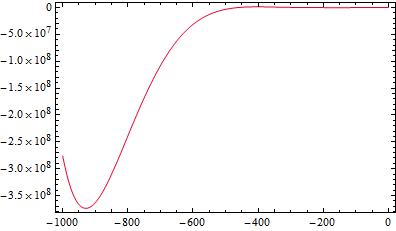}}
\caption{The figure of the solution on the interval
[-1000,0]}
    \end{center}
\end{figure}

\begin{figure}[H]
    \begin{center}
\vspace{0.3cm}
\scalebox{0.8}[0.55]{\includegraphics[trim=0mm 0mm
0mm 0mm, clip, width=9cm, height=12cm]{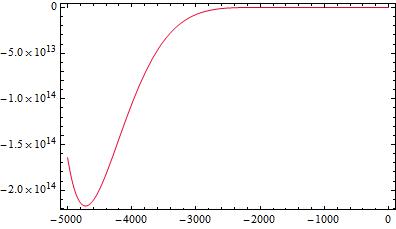}}
\caption{The figure of the solution on the interval
[-5000,0]}
    \end{center}
\end{figure}
On the other hand, on the interval $[0,\infty)$, the
solution is strictly increasing and goes to infinity as $t$ tends to
infinity.

For more interesting properties of the solution of the Cauchy problem (2.8), see next section.
\end{Theorem}
\end{Example}
\begin{Remark}
(\ref{2.10}) and Theorem 2.1 give the essential differences between (\ref{2.11}) and the traditional linear equation in (\ref{2.9}).
\end{Remark}

In fact, more generally, we consider the Cauchy problem for a system of functional
differential equations
\begin{equation}\label{2.6}
\left\{\aligned
&\aligned &\frac
{dx\left(t\right)}{dt}=A\left(t\right)x\left(t\right)+F\left
(x\left(\phi\left(t\right)\right)\right)x\left(t\right)+G\left(x\left
(\psi\left(t\right)\right)\right),\endaligned
\\
&x\left(0\right)=x_{0},\endaligned \right.\end{equation}
where $x(t)=(x_1(t),\cdots,x_n(t))^T$ is the unknown vector-valued function, $x_0$ is a constant vector. We have (see \cite{k-li})

\begin{Theorem} Let $A\left(t\right)$ be a continuous $ n\times n$
matrix-valued function of $t\in\mathbb{R}^{+}$, $F,\,G$ be locally
Lipschitzian continuous vector-valued functions of $x\in\mathbb{R}^n$, and
$\phi ,\,\psi$ be continuous real-valued functions of $ t\in\mathbb{R}^{+}$
with
$$0<\phi\left(t\right),\psi\left(t\right)<t\quad\text{ for all }\,\,
t>0.$$ Then Cauchy problem (\ref{2.6}) admits a unique global $C^1$
solution $ x=x\left(t\right)$ on $t\ge 0$.
\end{Theorem}

\begin{Example} We now consider the following Cauchy problem for the Burgers equation
\begin{equation}\label{2.7}
\left\{\aligned
&\aligned & u_t+\left(\frac{u^2}{2}\right)_x=0,\endaligned
\\
&t=0:\;\; u=\sin x. \endaligned \right.\end{equation}
It is well-known that the smooth solution of the Cauchy problem (\ref{2.7}) only exists on the strip $[0,1)\times \mathbb{R}$, and singularities, i.e., shock waves will appear at the time $t=1$ (e.g., at the point (1,0)). However, the Cauchy problem for the Burgers-like equation (\ref{1.7})
\begin{equation}\label{2.8}
\left\{\aligned
&\aligned & u_t(t,x)+u(t/2,x)u_x(t,x)=0,\endaligned
\\
&t=0:\;\; u=\sin x. \endaligned \right.\end{equation}
always admits a global smooth solution on the domain $\mathbb{R}^+\times\mathbb{R}$. See \cite{k-liu}. In fact, in the work \cite{k-liu}, the authors generalize the above results to the case of quasilinear hyperbolic systems of partial differential
equations.
\end{Example}

Clearly, Examples 2.1-2.3 give new philosophy enjoyed by functional differential equations (FDEs) and the essential differences between FDEs and traditional ones.

\section{Functional ordinary differential equation: Main results}

In this section, we investigate solutions to functional differential equations mentioned above.
For simplicity, in this section we only consider the solution of the following Cauchy problem
\begin{equation}\left\{\aligned
&\frac {dy\left(x\right)}{dx}=y\left(\frac x2\right),\\
&y\left(0\right)=1.\endaligned \right.\end{equation}
We shall prove the global existence of smooth solution, show regularities of zeros distribution and derive an exact estimate of the oscillation amplitude. In fact, our results can be easily generalized to other functional differential equations with a similar form.

Throughout this section, we denote the solution of the Cauchy problem (3.1) by $y=f(x)$. Direct calculation shows that the solution can be represented in power series form
\begin{equation}f(x) = \sum\limits_{n = 0}^\infty  {\frac{{{{x}^n}}}{{n!{2^{\frac{{n(n - 1)}}{2}}}}}}.
\end{equation}
In order to investigate the properties of the solution, we first prove the following useful lemma on the truncated polynomials of $f(x)$.
\begin{Lemma}For any integer $N > 2$, the truncated polynomial of $f(x)$, i.e.,\begin{equation}{P_N}(x) = \sum\limits_{n = 0}^N {\frac{{{x^n}}}{{n!{2^{\frac{{n(n - 1)}}{2}}}}}}
\end{equation}
has N-2 real zeros denoted by \begin{equation}{\{ {r_n}\} _{1 \leqslant n \leqslant N - 2}}\quad {\rm with} \quad  {r_{n + 1}} < {r_n} < 0,\end{equation} and a pair of imaginary zeros ${z_N},\overline {{z_N}} $. Moreover, the following estimates hold
\begin{equation}|{r_n}| < (n + 2){2^n}\quad (1 \le n \le N-2)
\end{equation}
and
\begin{equation}|{z_N}{|^2} > {2^N}.
\end{equation}
\end{Lemma}

By Lemma 3.1, we shall prove the following theorem on regularities of zeros distribution of the solution in Section 4.
\begin{Theorem}  Let $\mathbb{A}$ be the set of zeros of $f(z)$ in $\mathbb{C}$ .Then $\mathbb{A}\subset \mathbb{R}$ is countably infinite, and it can be represented as ${\{ {x_n}\} _{n \geq 1}}$ with the property $${x_{n + 1}} < {x_n} < 0.$$ Moreover, for any fixed integer $n \geq 1$, it holds that \begin{equation}{x_n} =  - (n + {\theta _n}){2^{n - 1}}\quad (0< \theta_n <1).
\end{equation}
\end{Theorem}

In fact, in next section we can show the following more precise estimate on the parameter ${\theta _n}$ in (3.7).
\begin{Remark}
For large $n$, we have
$$\frac{66}{25n} < {\theta _n} < \frac{2167}{789n}.$$
Moreover, it holds that
\begin{equation}\lim_{n \to \infty}n{\theta _n} = C\approx 2.744.
\end{equation}
\end{Remark}

Moreover, based on Theorem 3.1, we will prove the following theorem on the oscillation amplitude of the solution.
\begin{Theorem} Let $\mathbb{B}$ be the set of critical points of $y=f(x)$ in $\mathbb{R}$. Then $$\mathbb{B} = \{ 2x_1,2x_2,\cdots, 2x_n,\cdots\},$$
the extrema of $y=f(x)$ can be estimated as
\begin{equation}{\log _2}|f(2{x_n})| < \frac{1}{2}{n^2} + 2n,
\end{equation}
and their signs are alternate, i.e.,
\begin{equation}{( - 1)^n}f(2{x_n}) > 0.
\end{equation}
Moreover, for $x<0$ with large $|x|$, the following oscillation amplitude estimate holds
\begin{equation}{\log _2}|f(x)| < \frac{1}
{2}{\left( {{{\log }_2}|x|} \right)^2} + 2{\log _2}|x|.
\end{equation}
\end{Theorem}

\begin{Remark} On the other hand, by Gronwall inequality, we can derive the following rough estimate on oscillation amplitude for $x < 0$
\begin{equation}
{\log _2}|f(x)| < |x|\log _2e.
\end{equation}
Since $$|x|\gg {\left( {{{\log }_2}|x|} \right)^2}$$ for large $|x|$, the estimate (3.11) is much better than (3.12).
\end{Remark}

\section{Functional ordinary differential equation: Proof of Theorems 3.1 and 3.2}

In this section, we prove Theorems 3.1 and 3.2. As before, we still denote the solution by $y=f(x)$ throughout this section.

Recall that $f(x)$ in power series form reads
\begin{equation}f(x) = \sum\limits_{n = 0}^\infty  {\frac{{{{x}^n}}}{{n!{2^{\frac{{n(n - 1)}}{2}}}}}}.
\end{equation}
\begin{Proposition} For integer $k \ge 1$, the signs of values $f(-k{2^{k - 1}})$ and $f( - (k+1){2^{k - 1}})$ both are alternate, that is,
\begin{equation}{( - 1)^k}f( - k{2^{k - 1}}) < 0
\end{equation}
and
\begin{equation}{( - 1)^k}f(-(k+1){2^{k - 1}}) > 0.
\end{equation}
\end{Proposition}

\noindent{\bf Proof.} We first consider the value of the solution $y=f(x)$ at the point $-k2^{k - 1}$, namely, the value $f(-k{2^{k - 1}})$.

Direct calculation leads to
\begin{equation}f( - k{2^{k - 1}}) = \sum\limits_{n = 0}^\infty  {\frac{{{k^n}}}{{n!}}} {2^{nk - n(n + 1)/2}}{( - 1)^n}\mathop  = \limits^\Delta  \sum\limits_{n = 0}^\infty  {{( - 1)}^n{u_n}}.
\end{equation}
Then,
\begin{equation}\frac{{{u_{n + 1}}}}{{{u_n}}} = \frac{k}{{n + 1}}{2^{k - n - 1}}.
\end{equation}
Therefore,
\begin{equation}{u_n} {\rm \;\; is\;\; increasing\;\; for\;\; }n<k;\quad {\rm while}\;\; {u_n}  {\rm \;\; is\;\; decreasing\;\; for\;\;}  n>k.
\end{equation}
It is easy to check that
\begin{equation}{u_{2k}} = \frac{{{k^{2k}}}}{{(2k)!}}{2^{ - k}} < \frac{{{k^{2k}}}}{{{{\left( {\frac{{2k}}{e}} \right)}^{2k}}}}{2^{ - k}} = {\left( {\frac{e}{{2\sqrt 2 }}} \right)^{2k}} < 1.
\end{equation}
Based on this, we denote
\begin{equation}{v_j} = {u_{2k - j - 1}} - {u_j}\quad (0 \le j \le k - 1).
\end{equation}
We have the following lemma which will be proved later.
\begin{Lemma}The sequence ${\{ {v_j}\} _{0 \le j \le k - 1}}$ enjoys following properties
\begin{equation}{v_{k - 1}} = 0,\quad {v_j} > 0 \quad (0 \le j \le k - 2),
\end{equation}
\begin{equation}{v_j} < {v_{j + 1}} \quad (0 \le j \le k - 4)
\end{equation}
and
\begin{equation}{v_{k - 2}} - {v_{k - 3}} + {v_{k - 4}} - {v_{k - 5}} > 0.
\end{equation}
\end{Lemma}

On one hand, by Lemma 4.1 we obtain
\begin{equation}{( - 1)^k}\sum\limits_{n = 0}^{2k - 1} {{u_n}{{( - 1)}^n}}  = ( - {v_{k - 2}} + {v_{k - 3}} - {v_{k - 4}} + {v_{k - 5}}) - \sum\limits_{j = 0}^{k - 6} {{v_j}{{( - 1)}^{k - j}}}  < 0.
\end{equation}
On the other hand, direct calculation gives
\begin{equation}{v_0} = \frac{{{k^{2k - 1}}}}{{(2k - 1)!}} - 1 \ge 1\quad (k \ge 4)
\end{equation}
and
\begin{equation}{v_1} - {v_0} = (\frac{{{k^{2k - 2}}}}{{(2k - 2)!}} - k){2^{k - 1}} - \frac{{{k^{2k - 1}}}}{{(2k - 1)!}} + 1 \ge 1 \quad (k \ge 5).
\end{equation}
Hence,
\begin{equation}{( - 1)^k}\sum\limits_{n = 0}^{2k - 1} {{u_n}{{( - 1)}^n}}  <  - 1\quad (k \ge 4).
\end{equation}

Noting
\begin{equation}\left| {\sum\limits_{n = 2k}^\infty  {{u_n}{{( - 1)}^n}} } \right| < {u_{2k}} < 1,
\end{equation}
for $k\ge 4$ we have
\begin{equation}{( - 1)^k}\sum\limits_{n = 0}^\infty  {u{}_n{{( - 1)}^n} < } {( - 1)^k}\sum\limits_{n = 0}^{2k - 1} {u{}_n{{( - 1)}^n} + \left| {\sum\limits_{n = 2k}^\infty  {u{}_n{{( - 1)}^n}} } \right| < } {( - 1)^k}\sum\limits_{n = 0}^{2k - 1} {u{}_n{{( - 1)}^n} + 1 < } 0.
\end{equation}

Similarly, it is easy to check that the desired inequalities hold for $k=1,2,3$.
This finishes the proof of the first part of Proposition 4.1.

We next show the second part of Proposition 4.1. In fact, the method of proof is the same as above, here we only state some essential differences.

Direct calculation yields
\begin{equation}f( - (k + 1){2^{k - 1}}) = \sum\limits_{n = 0}^\infty  {\frac{{{{(k + 1)}^n}}}{{n!}}} {2^{nk - n(n + 1)/2}}{( - 1)^n}\mathop  = \limits^\Delta  \sum\limits_{n = 0}^\infty  {{u_n}{{( - 1)}^n}} .
\end{equation}
Then,
\begin{equation}\frac{{{u_{n + 1}}}}{{{u_n}}} = \frac{{k + 1}}{{n + 1}}{2^{k - n - 1}}.
\end{equation}
Therefore,
\begin{equation}
{u_n} {\rm \;\; is\;\; increasing\;\; for\;\; }n<k;\qquad {\rm while}\;\; {u_n}  {\rm \;\; is\;\; decreasing\;\; for\;\;}  n>k.
\end{equation}

Noting
\begin{equation}{u_2} = 1,\quad {u_{2k}} = \frac{{{{(k + 1)}^{2k}}}}{{(2k)!}}{2^{ - k}} < \frac{{{{(k + 1)}^{2k}}}}{{\sqrt {4\pi k} {{\left( {\frac{{2k}}{e}} \right)}^{2k}}}}{2^{ - k}} = \frac{1}{{\sqrt {4\pi k} }}{\left( {\frac{{e(k + 1)}}{{2\sqrt 2 k}}} \right)^{2k}} < 1\quad (k > 1)
\end{equation}
and denoting
\begin{equation}{v_j} = {u_{2k - j - 1}} - {u_j}\quad (0 \le j \le k - 1),
\end{equation}
in a manner similar to Lemma 4.1, we have the following lemma which will be proved later.
\begin{Lemma} The sequence ${\{ {v_j}\} _{0 \le j \le k - 1}}$ enjoys following properties
\begin{equation}{v_j} > 0\quad (0 \le j \le k - 1),
\end{equation}
\begin{equation}{v_j} < {v_{j + 1}}\quad (0 \le j \le k - 3)
\end{equation}
and
\begin{equation}{v_{k - 1}} - {v_{k - 2}} + {v_{k - 3}} - {v_{k - 4}} > 0.
\end{equation}
\end{Lemma}

By Lemma 4.2, similar to (4.12), we have
\begin{equation}{( - 1)^k}\sum\limits_{n = 0}^{2k - 1} {{u_n}{{( - 1)}^n}}  = ({v_{k - 1}} - {v_{k - 2}} + {v_{k - 3}} - {v_{k - 4}}) + \sum\limits_{j = 0}^{k - 5} {{v_j}{{( - 1)}^{k - j - 1}}}  > 0.
\end{equation}
Direct calculation yields
\begin{equation}{v_0} = \frac{{{{(k{\rm{ + 1}})}^{2k - 1}}}}{{(2k - 1)!}} - 1 \ge 1 \quad (k \ge {\rm{1}})
\end{equation}
and
\begin{equation}{v_1} - {v_0} = \left(\frac{{{{(k + 1)}^{2k - 2}}}}{{(2k - 2)!}} - k - 1\right){2^{k - 1}} - \frac{{{{(k + 1)}^{2k - 1}}}}{{(2k - 1)!}} + 1 \ge 1 (k \ge 4).
\end{equation}
Hence, it holds that
\begin{equation}{( - 1)^k}\sum\limits_{n = 0}^{2k - 1} {{u_n}{{( - 1)}^n}}  > 1 \quad (k \ge 4).
\end{equation}
Noting
\begin{equation}\left| {\sum\limits_{n = 2k}^\infty  {{u_n}{{( - 1)}^n}} } \right| < {u_{2k}} \le 1,
\end{equation}
we obtain
\begin{equation}{( - 1)^k}\sum\limits_{n = 0}^\infty  {{u_n}{{( - 1)}^n}}  > {( - 1)^k}\sum\limits_{n = 0}^{2k - 1} {{u_n}{{( - 1)}^n} - \left| {\sum\limits_{n = 2k}^\infty  {{u_n}{{( - 1)}^n}} } \right| > } {( - 1)^k}\sum\limits_{n = 0}^{2k - 1} {{u_n}{{( - 1)}^n} - 1 > } 0 \quad (k \ge 4).
\end{equation}

Moreover, it is easy to check that the inequality also holds for $k=1,2,3$. This proves the desired (4.3). Thus, Proposition 4.1 has been proved. $\quad\quad\quad\blacksquare$

\vskip 3mm

In what follows, we shall show Lemmas 4.1 and 4.2.

\noindent {\bf Proof of Lemma 4.1.}  On one hand, we notice
\begin{equation}{v_{k - 1}} = {u_k} - {u_{k - 1}} = 0.
\end{equation}
On the other hand, we have
\begin{equation}\prod\limits_{i = 1}^{2k - 1 - 2j} {(j + i)}  < {k^{2k - 1 - 2j}}\quad (0 \le j \le k - 2).
\end{equation}
Hence,
\begin{equation}\frac{{{k^j}}}{{j!}} < \frac{{{k^{2k - 1 - j}}}}{{(2k - 1 - j)!}}\quad (0 \le j \le k - 2),
\end{equation}
namely,
\begin{equation}{\rm{ }}{u_j} < {u_{2k - 1 - j}}\quad (0 \le j \le k - 2)
\end{equation}
and
\begin{equation}{v_j} > 0 \quad (0 \le j \le k - 2).
\end{equation}
This proves (4.9).

In order to prove (4.10), namely,
\begin{equation}{v_j} < {v_{j + 1}}\quad (0 \le j \le k - 4),
\end{equation}
by direct calculation it suffices to show
\begin{equation}\left({2^{k - j - 1}} - \frac{{j + 1}}{k}\right){w_j} < {2^{k - j - 1}} - \frac{k}{{2k - 1 - j}} \quad (0 \le j \le k - 4),
\end{equation}
where
\begin{equation}{w_j} =  \frac{{(2k - 2 - j)!}}{{(j + 1)!{k^{2k - 3 - 2j}}}}.
\end{equation}

Note that
\begin{equation}\frac{{{w_{j + 1}}}}{{{w_j}}} = \frac{{{k^2}}}{{(j + 2)(2k - 2 - j)}} > 1.
\end{equation}
This is to say, $w_j$ is increasing. Hence, we only need to prove
\begin{equation}\left({2^{k - j - 1}} - \frac{{j + 1}}{k}\right){w_{k - 4}} < {2^{k - j - 1}} - \frac{k}{{2k - 1 - j}}\quad (0 \le j \le k - 4),
\end{equation}
i.e.,
\begin{equation}\left({2^{k - j - 1}} - \frac{{j + 1}}{k}\right)\frac{{({k^2} - 1)({k^2} - 4)}}{{{k^4}}} < {2^{k - j - 1}} - \frac{k}{{2k - 1 - j}}\quad (0 \le j \le k - 4)
\end{equation}
and
\begin{equation}\frac{1}{{2k - 1 - j}} - \frac{{(j + 1)({k^2} - 1)({k^2} - 4)}}{{{k^6}}} - {2^{k - j - 1}}\frac{{5{k^2} - 4}}{{{k^5}}} < 0\quad (0 \le j \le k - 4).
\end{equation}

Denote
\begin{equation}t = 2k - 1 - j\quad (k + 3 \le t \le 2k - 1)
\end{equation}
and
\begin{equation} g(t) = \frac{1}{t} - \frac{{(2k - t)({k^2} - 1)({k^2} - 4)}}{{{k^6}}} - {2^{t - k}}\frac{{5{k^2} - 4}}{{{k^5}}}.
\end{equation}

In what follows, we shall prove $$g(t)< 0,$$ which will leads to our desired result immediately.

In fact, direct calculation yields
\begin{equation}g'(t) =  - \frac{1}{{{t^2}}} + \frac{{({k^2} - 1)({k^2} - 4)}}{{{k^6}}} - {2^{t - k}}\frac{{5{k^2} - 4}}{{{k^5}}}\ln 2
\end{equation}
and
\begin{equation}g''(t) = \frac{2}{{{t^3}}} - {2^{t - k}}\frac{{5{k^2} - 4}}{{{k^5}}}{(\ln 2)^2}.
\end{equation}
Combining
\begin{equation}g''(k + 3) = \frac{2}{{{{(k + 3)}^3}}} - \frac{{8(5{k^2} - 4)}}{{{k^5}}}{(\ln 2)^2} =  - \frac{{40{{(\ln 2)}^2} - 2}}{{{k^3}}} + O(k^{-4}) < 0,
\end{equation}
\begin{equation}g'(k + 3) =  - \frac{1}{{{{(k + 3)}^2}}} + \frac{{({k^2} - 1)({k^2} - 4)}}{{{k^6}}} - \frac{{8(5{k^2} - 4)}}{{{k^5}}}\ln 2 =  - \frac{{40\ln 2 - 6}}{{{k^3}}} + O(k^{-4}) < 0
\end{equation}
and
\begin{equation}g(k + 3) = \frac{1}{{k + 3}} - \frac{{(k - 3)({k^2} - 1)({k^2} - 4)}}{{{k^6}}} - \frac{{8(5{k^2} - 4)}}{{{k^5}}} =  - \frac{{26}}{{{k^3}}} + O(k^{-4}) < 0
\end{equation}
gives
\begin{equation}g(t) \le g(k + 3) < 0.
\end{equation}
This gives the desired (4.10).

Finally, we prove (4.11).

Noting
\begin{equation}
 {v_{k - 3}} - {v_{k - 2}} = {u_{k + 2}} - {u_{k - 3}} - {u_{k + 1}} + {u_{k - 2}}
  = \frac{{{k^{k{\rm{ - }}4}}}}{{(k + 4)!}}\left( {128{k^6} + O({k^5})} \right){2^{\frac{{{k^2} - k - 20}}{2}}}
\end{equation}
and
\begin{equation}
 {v_{k - 4}} - {v_{k - 5}} = {u_{k + 3}} - {u_{k - 4}} - {u_{k + 4}} + {u_{k - 5}}
  = \frac{{{k^{k - 4}}}}{{(k + 4)!}}\left( {194{k^6} + O({k^5})} \right){2^{\frac{{{k^2} - k - 20}}{2}}},
\end{equation}
we obtain
\begin{equation}{v_{k - 2}} - {v_{k - 3}} + {v_{k - 4}} - {v_{k - 5}} = \frac{{{k^{k - 4}}}}{{(k + 4)!}}\left( {66{k^6} + O({k^5})} \right){2^{\frac{{{k^2} - k - 20}}{2}}} > 0,
\end{equation}
which is nothing but the desired estimate (4.11). Thus, the proof of Lemma 4.1 is completed. $\quad\quad\quad\blacksquare$

\vskip 4mm

We next prove Lemma 4.2.

\noindent {\bf Proof of Lemma 4.2.}
Noting
\begin{equation}\prod\limits_{i = 1}^{2k - 1 - 2j} {(j + i)}  < {(k + 1)^{2k - 1 - 2j}}\quad (0 \le j \le k - 1),
\end{equation}
we have
\begin{equation}\frac{{{{(k + 1)}^j}}}{{j!}} < \frac{{{{(k + 1)}^{2k - 1 - j}}}}{{(2k - 1 - j)!}}\quad (0 \le j \le k - 1).
\end{equation}
This implies that
\begin{equation}{u_j} < {u_{2k - 1 - j}}\quad (0 \le j \le k - 1),
\end{equation}
and then
\begin{equation}{v_j} > 0\quad (0 \le j \le k - 1).
\end{equation}
This is nothing but the desired (4.23).

In order to prove (4.24), namely,
\begin{equation}{v_j} < {v_{j + 1}}\quad (0 \le j \le k - 3),
\end{equation}
by direct calculations it suffices to show
\begin{equation}\left({2^{k - j - 1}} - \frac{{j + 1}}{{k + 1}}\right){w_j} < {2^{k - j - 1}} - \frac{{k + 1}}{{2k - 1 - j}}\quad (0 \le j \le k - 3),
\end{equation}
where $w_j$ is given by
\begin{equation}{w_j}=\frac{{(2k - 2 - j)!}}{{(j + 1)!{{(k + 1)}^{2k - 3 - 2j}}}}.
\end{equation}
Moreover, noting that
\begin{equation}\frac{{{w_{j + 1}}}}{{{w_j}}} = \frac{{{{(k + 1)}^2}}}{{(j + 2)(2k - 2 - j)}} > 1,
\end{equation}
that is, $w_j$ is increasing, we only need to prove
\begin{equation}\left({2^{k - j - 1}} - \frac{{j + 1}}{{k + 1}}\right){w_{k - 3}} < {2^{k - j - 1}} - \frac{{k + 1}}{{2k - 1 - j}}\quad (0 \le j \le k - 3),
\end{equation}
namely,
\begin{equation}\left({2^{k - j - 1}} - \frac{{j + 1}}{{k + 1}}\right)\frac{{k(k - 1)}}{{{{(k + 1)}^2}}} < {2^{k - j - 1}} - \frac{{k + 1}}{{2k - 1 - j}}\quad (0 \le j \le k - 3)
\end{equation}
or
\begin{equation}\frac{1}{{2k - 1 - j}} - \frac{{(j + 1)k(k - 1)}}{{{{(k + 1)}^4}}} - {2^{k - j - 1}}\frac{{3k + 1}}{{{{(k + 1)}^3}}} < 0\quad (0 \le j \le k - 3).
\end{equation}

Denote
\begin{equation} t \mathop  =  2k - 1 - j\quad (k + 2 \le t \le 2k - 1),
\end{equation}
\begin{equation}g(t) \mathop  =  \frac{1}{t} - \frac{{(2k - t)k(k - 1)}}{{{{(k + 1)}^4}}} - {2^{t - k}}\frac{{3k + 1}}{{{{(k + 1)}^3}}}.
\end{equation}
We next prove $$g(t)< 0,$$ which leads to our desired result immediately.

In fact, direct calculations show that
\begin{equation}g'(t) =  - \frac{1}{{{t^2}}} + \frac{{k(k - 1)}}{{{{(k + 1)}^4}}} - {2^{t - k}}\frac{{3k + 1}}{{{{(k + 1)}^3}}}\ln 2
\end{equation}
and
\begin{equation}g''(t) = \frac{2}{{{t^3}}} - {2^{t - k}}\frac{{3k + 1}}{{{{(k + 1)}^3}}}{(\ln 2)^2}.
\end{equation}
On the other hand, noting
\begin{equation}g''(k + 2) = \frac{2}{{{{(k + 2)}^3}}} - \frac{{4(3k + 1)}}{{{{(k + 1)}^3}}}{(\ln 2)^2} =  - \frac{{12{{(\ln 2)}^2}}}{{{k^2}}} + O(k^{-3})<0,
\end{equation}
\begin{equation}g'(k + 2) =  - \frac{1}{{{{(k + 2)}^2}}} + \frac{{k(k - 1)}}{{{{(k + 1)}^4}}} - \frac{{4(3k + 1)}}{{{{(k + 1)}^3}}}\ln 2 =  - \frac{{12\ln 2}}{{{k^2}}} + O(k^{-3}) < 0,
\end{equation}
and using
\begin{equation}g(k + 2) = \frac{1}{{k + 2}} - \frac{{k(k - 1)(k - 2)}}{{{{(k + 1)}^4}}} - \frac{{4(3k + 1)}}{{{{(k + 1)}^3}}} =  - \frac{7}{{{k^2}}} + O(k^{-3}) < 0,
\end{equation}
we obtain
\begin{equation}g(t) \le g(k + 2) < 0.
\end{equation}
This proves (4.24).

Finally, we show (4.25).

In fact, noting
\begin{equation}
 {v_{k - 2}} - {v_{k - 1}} = {u_{k + 1}} - {u_{k - 2}} - {u_k} + {u_{k - 1}}= \frac{{{{(k + 1)}^{k - 4}}}}{{(k + 3)!}}\left(32{k^6} + O({k^5})\right){2^{\frac{{{k^2} - k - 12}}{2}}}
\end{equation}
and using
\begin{equation}
 {v_{k - 3}} - {v_{k - 4}} = {u_{k + 2}} - {u_{k - 3}} - {u_{k + 3}} + {u_{k - 4}}= \frac{{{{(k + 1)}^{k - 4}}}}{{(k + 3)!}}\left(33{k^6} + O({k^5})\right){2^{\frac{{{k^2} - k - 12}}{2}}},
\end{equation}
we have
\begin{equation}{v_{k - 1}} - {v_{k - 2}} + {v_{k - 3}} - {v_{k - 4}} = \frac{{{{(k + 1)}^{k - 4}}}}{{(k + 3)!}}\left( {{k^6} + O({k^5})} \right){2^{\frac{{{k^2} - k - 12}}{2}}} > 0,
\end{equation}
which leads to the desired (4.25) directly. Thus, the proof of Lemma 4.2 is finished.  $\quad\quad\quad\blacksquare$

\vskip 4mm
In a similar way, we can show the following lemma on more precise results on the sign of function values.
\begin{Lemma} For large $k$, it holds that
\begin{equation}{( - 1)^k}f\left( - \left(k + \frac{\lambda }{k}\right){2^{k - 1}}\right) < 0,\quad \lambda  \in \left(0,\frac{{66}}{{25}}\right)
\end{equation}
and
\begin{equation}{( - 1)^k}f\left( - \left(k + \frac{\lambda }{k}\right){2^{k - 1}}\right) > 0,\quad \lambda  \in \left(\frac{{2167}}{{789}},\infty \right).
\end{equation}
\end{Lemma}
\noindent{\bf Proof.} Since the proof is almost the same as above, we only prove it in a brief way.

For simplicity, we assume $\lambda  > 1$. By direct calculations, we have
\begin{equation}f\left( - \left(k + \frac{\lambda }{k}\right){2^{k - 1}}\right) = \sum\limits_{n = 0}^\infty {{( - 1)}^n} {\frac{{{{(k + \frac{\lambda }{k})}^n}}}{{n!}}{2^{nk - n(n + 1)/2}}} \mathop  = \limits^\Delta  \sum\limits_{n = 0}^\infty {{( - 1)}^n} {u_n}
\end{equation}
Similarly, we have
\begin{equation}
{u_n} {\rm \;\; is\;\; increasing\;\; for\;\; }n<k;\quad {u_n}  {\rm \;\; is\;\; decreasing\;\; for\;\;}  n>k; \quad u_{2k}<1.
\end{equation}
Denote
\begin{equation}{v_j}  \mathop  = {u_{2k - 1 - j}} - {u_j}\quad (0 \le j \le k - 1).
\end{equation}
By similar calculation and estimation, we get
\begin{equation}{v_j} > 0\quad (0 \le j \le k - 1)
\end{equation}
and
\begin{equation}{v_j} < {v_{j + 1}}\quad (0 \le j \le k - 3).
\end{equation}
It is easy to see that, for  large $k$
\begin{equation}{v_{k - 1}} - {v_{k - 2}} + {v_{k - 3}} - {v_{k - 4}} + {v_{k - 5}} = \frac{{{{(k + \frac{\lambda }{k})}^{k + 2}}}}{{(k + 4)!}}\left( {25\lambda  - 66 + O\left( k^{-1} \right)} \right){2^{\frac{{{k^2} - k - 20}}{2}}}
\end{equation}
and
\begin{equation}{v_{k - 1}} - {v_{k - 2}} + {v_{k - 3}} - {v_{k - 4}} + {v_{k - 5}} - {v_{k - 6}} = \frac{{{{(k + \frac{\lambda }{k})}^{k + 3}}}}{{(k + 5)!}}\left( {789\lambda  - 2167 + O\left( k^{-1} \right)} \right){2^{\frac{{{k^2} - k - 30}}{2}}}.
\end{equation}
Hence, for large $k$ it holds that
\begin{equation}{( - 1)^k}\sum\limits_{n = 0}^{2k - 1} {{( - 1)}^n}{{u_n}}  = ({v_{k - 1}} - {v_{k - 2}} + {v_{k - 3}} - {v_{k - 4}} + {v_{k - 5}}) - \sum\limits_{j = 0}^{k - 6} {{{( - 1)}^{k - j}}{v_j}}  < 0\quad {\rm for }\;\;\lambda  < \frac{{66}}{{25}}
\end{equation}
and
\begin{equation}{( - 1)^k}\sum\limits_{n = 0}^{2k - 1} {{( - 1)}^n{u_n}{}}  = ({v_{k - 1}} - {v_{k - 2}} + {v_{k - 3}} - {v_{k - 4}} + {v_{k - 5}} - {v_{k - 6}}) + \sum\limits_{j = 0}^{k - 7} {{( - 1)}^{k - j - 1}{v_j}{}}  > 0\quad {\rm for }\;\;\lambda  > \frac{{2167}}{{789}}.
\end{equation}

Similarly, we can show
\begin{equation}\left| {\sum\limits_{n = 0}^{2k - 1} {{( - 1)}^n{u_n}{}} } \right| > 1
\end{equation}
and
\begin{equation}\left| {\sum\limits_{n = 2k}^\infty  {{( - 1)}^n{u_n}} } \right| < 1.
\end{equation}
This proves Lemma 4.3. $\quad\quad\quad\blacksquare$

\vskip 4mm
Again, notice
\begin{equation}{v_{k - j}} = \frac{{{{(k + \frac{\lambda }
{k})}^{k + j - 3}}}}
{{(k + j - 1)!}}\left[ {(2j - 1)\lambda  + \frac{j}
{6}(j - 1)(2j - 1) + O\left( k^{-1} \right)} \right]{2^{\frac{{{k^2} - k - j(j - 1)}}
{2}}}.
\end{equation}
Fixing any integer $m$ with $2m<k$, when $k$ is large enough, we have
\begin{equation}\begin{gathered}
  \sum\limits_{j = 1}^{2m - 1} {{( - 1)}^{j - 1}{v_{k - j}}}  =  \hfill \\
  \frac{{{{(k + \frac{\lambda }
{k})}^{k - 2}}}}
{{k!}}\left[ {\lambda \sum\limits_{j = 1}^{2m - 1} {(2j - 1){( - 1)}^{j - 1}{2^{\frac{{j - {j^2}}}
{2}}}}  + \sum\limits_{j = 1}^{2m - 1} {\frac{j}
{6}(j - 1)(2j - 1){( - 1)}^{j - 1}{2^{\frac{{j - {j^2}}}
{2}}}}  + O\left( k^{-1} \right)} \right]{2^{\frac{{{k^2} - k}}
{2}}} \hfill \\
\end{gathered}
\end{equation}
and
\begin{equation}\begin{gathered}
  \sum\limits_{j = 1}^{2m} {{( - 1)}^{j - 1}{v_{k - j}}}  =  \hfill \\
  \frac{{{{(k + \frac{\lambda }
{k})}^{k - 2}}}}
{{k!}}\left[ {\lambda \sum\limits_{j = 1}^{2m} {(2j - 1){2^{\frac{{j - {j^2}}}
{2}}}{{( - 1)}^{j - 1}}}  + \sum\limits_{j = 1}^{2m} {\frac{j}
{6}(j - 1)(2j - 1){2^{\frac{{j - {j^2}}}
{2}}}{{( - 1)}^{j - 1}}}  + O\left( k^{-1} \right)} \right]{2^{\frac{{{k^2} - k}}
{2}}}. \hfill \\
\end{gathered}
\end{equation}
Hence
\begin{equation}{( - 1)^k}\sum\limits_{n = 0}^{2k - 1} {{( - 1)}^n{u_n}}  < 0, \quad \lambda  < \frac{{\sum\limits_{j = 1}^{2m - 1} {\frac{j}
{6}(j - 1)(2j - 1){2^{\frac{{j - {j^2}}}
{2}}}{{( - 1)}^j}} }}
{{\sum\limits_{j = 1}^{2m - 1} {(2j - 1){2^{\frac{{j - {j^2}}}
{2}}}{{( - 1)}^{j - 1}}} }}
\end{equation}
and
\begin{equation}{( - 1)^k}\sum\limits_{n = 0}^{2k - 1} {{( - 1)}^n{u_n}}  > 0,\quad \lambda  > \frac{{\sum\limits_{j = 1}^{2m} {\frac{j}
{6}(j - 1)(2j - 1){2^{\frac{{j - {j^2}}}
{2}}}{{( - 1)}^j}} }}
{{\sum\limits_{j = 1}^{2m} {(2j - 1){2^{\frac{{j - {j^2}}}
{2}}}{{( - 1)}^{j - 1}}} }}.
\end{equation}
Moreover, we have
\begin{equation}\mathop {\lim }\limits_{m \to \infty } \frac{{\sum\limits_{j = 1}^{2m} {\frac{j}
{6}(j - 1)(2j - 1){2^{\frac{{j - {j^2}}}
{2}}}{{( - 1)}^j}} }}
{{\sum\limits_{j = 1}^{2m} {(2j - 1){2^{\frac{{j - {j^2}}}
{2}}}{{( - 1)}^{j - 1}}} }} = C \approx 2.744.
\end{equation}
This result will be used in the following discussion.

\vskip 4mm
By Proposition 4.1, we can prove
\begin{Proposition} The solution $f(x)$ has a countably infinite number of different real zeros denoted  by $\mathbb{A} \mathop  = \{x_n\} _{n \ge 1}$ with ${x_{n + 1}} < {x_n} < 0$. Moreover it holds that
\begin{equation}{x_n} =  - (n + {\theta _n}){2^{n - 1}}\quad (0 < {\theta _n} < 1).
\end{equation}
\end{Proposition}
\noindent{\bf Proof.}  Noting that the convergence radius of $f(z)$ in $\mathbb{C}$ is $+ \infty$ and $f(z)\not  \equiv 0$, we see that $f(z)$ is a nonzero analytic function. Since the zeros of a nonzero analytic function are isolated, we observe that $f(x)$ has an at most countably infinite number of real zeros.

By Proposition 4.1 and intermediate value theorem, $f(x)$ has an infinite number of real zeros. Hence we conclude that $f(x)$ has a countably infinite number of real zeros. Since the real zeros are isolated and negative, the set of real zeros of $f(x)$ can be uniquely represented as ${\{ {x_n}\} _{n \ge 1}}$ with ${x_{n + 1}} < {x_n} < 0$.

In what follows, we prove the estimate of ${x_n}$.

Notice that ${\{ 2{x_n}\} _{n \geq 1}}$ is the set of critical points of $f(x)$. Then for any integer $k \geq 1$, there are $k-1$ critical points in the interval $(2{x_k}, + \infty )$. By mean value theorem, we observe that $f(x)$ has at most $k$ different real zeros and when the number of zeros is $k$, every zero is simple. Again, by noting
\begin{equation}{\{ {x_n}\} _{1 \leqslant n \leqslant k}} \subset (2{x_k}, + \infty ),
\end{equation}
we observe immediately that $f(x)$ has exactly $k$ different real zeros in $(2{x_k}, + \infty )$ and every zero is simple. Hence, we get
\begin{equation}{x_{n + 1}} < 2{x_n}\quad  (n \geq 1).
\end{equation}

On one hand, for any fixed $x \in [ - 1,0)$, it holds that
\begin{equation}f(x) = \sum\limits_{n = 0}^\infty  {\frac{{{x^n}}}
{{n!{2^{\frac{{n(n - 1)}}
{2}}}}}}  = \sum\limits_{n = 0}^\infty  {\frac{{|x{|^n}}}
{{n!{2^{\frac{{n(n - 1)}}
{2}}}}}} {( - 1)^n}\mathop  = \limits^\Delta  \sum\limits_{n = 0}^\infty  {{a_n}{{( - 1)}^n}}.
\end{equation}
Since $\{a_n\}$ is decreasing and
\begin{equation}{a_0} > 0,
\end{equation}
we obtain
\begin{equation}f(x) > 0,\quad \forall\;x \in [-1,0).
\end{equation}
Therefore, $${x_1} <  - 1.$$
On the other hand, by Proposition 4.1, $f(x)$ has zeros in interval $( - 2, - 1)$. Hence
\begin{equation}{x_1} \in ( - 2, - 1).
\end{equation}

In what follows, we finish the proof by induction.

To do so, we assume that
\begin{equation}{x_k} \in ( - (k + 1){2^{k - 1}}, - k{2^{k - 1}})
\end{equation}holds for integer $k \geq 1$. Then
\begin{equation}{x_{k + 1}} < 2{x_k} <  - k{2^k} \leqslant  - (k + 1){2^{k - 1}}.
\end{equation}
We next prove \begin{equation}{x_{k + 1}} <  - (k + 1){2^k}
\end{equation} by contradiction.

In fact, assuming
\begin{equation}{x_{k + 1}} \in ( - (k + 1){2^k}, - (k + 1){2^{k - 1}}),
\end{equation}£¬
then by Proposition 4.1, we have
\begin{equation}f( - (k + 1){2^k})f( - (k + 1){2^{k - 1}}) > 0.
\end{equation}
Since ${x_{k + 1}}$ is simple, there exists a point
\begin{equation}{x^*} \in ( - (k + 1){2^k},{x_{k + 1}})
\end{equation}such that
\begin{equation}f({x^*})f( - (k + 1){2^k}) < 0.
\end{equation}
Therefore, $f(x)$ has zeros in $( - (k + 1){2^k},{x^*})$. Hence, we have
\begin{equation}{x_{k + 2}} \in ( - (k + 1){2^k},{x_{k + 1}}),
\end{equation}
which contradicts to
\begin{equation}{x_{k + 2}} < 2{x_{k + 1}} <  - (k + 1){2^k}.
\end{equation}
Consequently, we get
\begin{equation}{x_{k + 1}} <  - (k + 1){2^k}.
\end{equation}
Since $f(x)$ has zeros in the interval $( - (k + 2){2^k}, - (k + 1){2^k})$, we obtain
\begin{equation}{x_{k + 1}} \in ( - (k + 2){2^k}, - (k + 1){2^k}).
\end{equation}
Thus, the desired result has been established for any integer $k \geq 1$. This finishes the proof of Proposition 4.2. $\quad\quad\quad\blacksquare$

The next proposition is on the truncated polynomials.
\begin{Proposition} For any given integer $N > 2$, the following truncated polynomial of $f(x)$
\begin{equation}{P_N}(x) = \sum\limits_{n = 0}^N {\frac{{{x^n}}}{{n!{2^{\frac{{n(n - 1)}}{2}}}}}}
\end{equation}
has a pair of imaginary zeros $\{{z_N},\overline {{z_N}}\}$ and $N-2$ real zeros, which are denoted by
\begin{equation}{\{ {r_n}\} _{1 \leqslant n \leqslant N - 2}}\quad {\rm with} \quad \ {r_{n + 1}} < {r_n} < 0.
\end{equation}
Moreover, the following estimates hold
\begin{equation}|{r_n}| < (n + 2){2^n}\quad (1 \le n \le N - 2)
\end{equation}
and
\begin{equation}|{z_N}{|^2} > {2^N}.
\end{equation}
\end{Proposition}
\noindent{\bf Proof.} The proof consists of six steps. We denote ${\xi _k} =  - (k + 1){2^{k - 1}}(k \geq 0)$.

{\bf Step I.}$\;\;$ For $k \in \left\{1,\cdots,\left[ {\frac{{N + 1}}{2}} \right] \right\}$, we claim
\begin{equation}{( - 1)^k}{P_N}( {\xi _k}) > 0.
\end{equation}

In fact, by direct calculation, we have
\begin{equation}{P_N}( - (k + 1){2^{k - 1}}) = \sum\limits_{n = 0}^N  {{( - 1)}^n\frac{{{{(k + 1)}^n}}}{{n!}}{2^{nk - n(n + 1)/2}}} \mathop  = \limits^\Delta  \sum\limits_{n = 0}^N  {{( - 1)}^n{u_n}}.
\end{equation}
Then
\begin{equation}\frac{{{u_{n + 1}}}}{{{u_n}}} = \frac{{k + 1}}{{n + 1}}{2^{k - n - 1}}.
\end{equation}
Hence
\begin{equation}
{u_n} {\rm \;\; is\;\; increasing\;\; for\;\; }n<k;\quad {u_n}  {\rm \;\; is\;\; decreasing\;\; for\;\;}  n>k; \quad u_{2k}<1.
\end{equation}
As before, we denote
\begin{equation}{v_j} = {u_{2k - 1 - j}} - {u_j},
\end{equation}
and then we obtain
\begin{equation}{v_j} > 0\quad (0 \le j \le k - 1),
\end{equation}
\begin{equation}{v_j} < {v_{j + 1}}\quad (0 \le j \le k - 3)
\end{equation}
and
\begin{equation}{v_{k - 1}} - {v_{k - 2}} + {v_{k - 3}} - {v_{k - 4}} > 0.
\end{equation}
Again, we have
\begin{equation}{( - 1)^k}\sum\limits_{n = 0}^{2k - 1} {{( - 1)}^n{u_n}}  > 1.
\end{equation}
Since $2k - 1 \le N$ and
\begin{equation}\left| {\sum\limits_{n = 2k}^{N} {{( - 1)}^n{u_n}} } \right| < 1,
\end{equation}we get
\begin{equation}{( - 1)^k}\sum\limits_{n = 0}^N {{( - 1)}^n{u_n}}  > 0.
\end{equation}

{\bf Step II.}$\;\;$ For any given integer $k = N - 2j$ $(j \in {\mathbb{N}_ + })$ with $2k - 1 > N$,
we claim
\begin{equation}{( - 1)^N}{P_N}( {\xi _k}) > 0.
\end{equation}

Similar to Step I, denote
\begin{equation}{v_j} = {u_{2k - 1 - j}} - {u_j}.
\end{equation}
Then it holds that
\begin{equation}{v_j} > 0\quad (0 \le j \le k - 1),
\end{equation}
\begin{equation}{v_j} < {v_{j + 1}}\quad (0 \le j \le k - 3)
\end{equation}
and
\begin{equation}{v_{k - 1}} - {v_{k - 2}} + {v_{k - 3}} - {v_{k - 4}} > 0.
\end{equation}
Then
\begin{equation}\begin{gathered}
  {( - 1)^N}{P_N}( - (k + 1){2^{k - 1}}) = \sum\limits_{n = 2k - 1 - N}^{k - 1} {{v_n}{{( - 1)}^{n - k + 1}}}  + \sum\limits_{n = 0}^{2k - 2- N} {{u_n}{{( - 1)}^{N - n}}}  \hfill \\
   \qquad\qquad\qquad\qquad\qquad\quad \geq \sum\limits_{n = 2k - 1 - N}^{k - 1} {{v_n}{{( - 1)}^{n - k + 1}}}  > {v_{2k - 2 - N}} > 0. \hfill \\
\end{gathered}
\end{equation}
This finishes the proof of Step II.

{\bf Step III.}$\;\;$  For any integer $N>7$, we claim
\begin{equation}
{( - 1)^{N - 1}}{P_N}( - 5N{2^{N - 6}}) > 0.
\end{equation}

In fact, direct calculations yield
\begin{equation}{( - 1)^{N - 1}}{P_N}( - 5N{2^{N - 6}}) = \sum\limits_{n = 0}^N {( - 1)^{n - N + 1}}{\frac{{{{(5N{2^{N - 6}})}^n}}}
{{n!{2^{\frac{{n(n - 1)}}
{2}}}}}} \mathop  = \limits^\Delta  \sum\limits_{n = 0}^N {( - 1)^{n - N + 1}}{{a_n}}.
\end{equation}
Then
\begin{equation}\frac{{{a_n}}}
{{{a_{n - 1}}}} = \frac{{5N{2^{N - 6}}}}
{{(n + 1){2^n}}}.
\end{equation}
Hence
\begin{equation}
{a_n} {\rm \;\; is\;\; increasing\;\; for\;\; }n<N-3;\quad{\rm while}\;\; {a_n}  {\rm \;\; is\;\; decreasing\;\; for\;\;}  n>N-3.
\end{equation}
Notice
\begin{equation}\begin{gathered}
  {( - 1)^{N - 1}}\sum\limits_{n = 0}^N {( - 1)^n}{{a_n}}  \hfill \\
   =  - \left[ {({a_N} - {a_{N - 7}}) - ({a_{N - 1}} - {a_{N - 6}}) + ({a_{N - 2}} - {a_{N - 5}}) - ({a_{N - 3}} - {a_{N - 4}}) + \sum\limits_{n = N - 8}^N {{( - 1)}^{n - N}}{{a_n}} } \right], \hfill \\
\end{gathered}
\end{equation}
where
\begin{equation}{a_N} - {a_{N - 7}} = \frac{{{{(5N{2^{N - 6}})}^{N - 7}}}}
{{N!{2^{\frac{{(N - 7)(N - 8)}}
{2}}}}}\left[ {\frac{{61741}}
{{16384}}{N^7} + O({N^6})} \right],
\end{equation}
\begin{equation}{a_{N - 1}} - {a_{N - 6}} = \frac{{{{(5N{2^{N - 6}})}^{N - 6}}}}
{{N!{2^{\frac{{(N - 6)(N - 7)}}
{2}}}}}\left[ {\frac{{2101}}
{{1024}}{N^5} + O({N^4})} \right],
\end{equation}
\begin{equation}{a_{N - 2}} - {a_{N - 5}} = \frac{{{{(5N{2^{N - 6}})}^{N - 5}}}}
{{N!{2^{\frac{{(N - 5)(N - 6)}}
{2}}}}}\left[ {\frac{{61}}
{{64}}{N^3} + O({N^2})} \right]
\end{equation}
and
\begin{equation}{a_{N - 3}} - {a_{N - 4}} = \frac{{{{(5N{2^{N - 6}})}^{N - 3}}}}
{{N!{2^{\frac{{(N - 3)(N - 4)}}
{2}}}}}\left[ {\frac{1}
{4}N + O(1)} \right].
\end{equation}
Direct calculation leads to
\begin{equation}\begin{gathered}
  ({a_N} - {a_{N - 7}}) - ({a_{N - 1}} - {a_{N - 6}}) + ({a_{N - 2}} - {a_{N - 5}}) - ({a_{N - 3}} - {a_{N - 4}}) \hfill \\
  \qquad\qquad = \frac{{{{(5N{2^{N - 6}})}^{N - 7}}}}
{{N!{2^{\frac{{(N - 7)(N - 8)}}
{2}}}}}\left[ { - \frac{{5619}}
{{16384}}{N^7} + O({N^6})} \right] <  - \frac{{{{(5N{2^{N - 6}})}^{N - 8}}}}
{{(N - 8)!{2^{\frac{{(N - 8)(N - 9)}}
{2}}}}} =  - {a_{N - 8}}. \hfill \\
\end{gathered}
\end{equation}
Therefore, we obtain
\begin{equation}{( - 1)^{N - 1}}\sum\limits_{n = 0}^N {( - 1)^n}{{a_n}}  > 0.
\end{equation}
This proves (4.135).

{\bf Step IV.}$\;\;$  For any given integer $N>2$, we claim that ${P_N}(x)$ has imaginary zeros.
In what follows, we prove this statement by contradiction.

Assume that ${P_N}(x)$ only has $N$ real zeros. Then by mean-value theorem, we observe that
\begin{equation}\frac{{{d^{N - 3}}}}{{d{x^{N - 3}}}}{P_N}(x)
\end{equation}has $3$ real zeros.
 Direct calculation yields
\begin{equation}\frac{{{d^{N - 3}}}}
{{d{x^{N - 3}}}}P{}_N(x) = {2^{ - \frac{{(n - 3)(n - 4)}}
{{\text{2}}}}} + {2^{ - \frac{{(n - 2)(n - 3)}}
{{\text{2}}}}}x + \frac{1}
{2}{2^{ - \frac{{(n - 1)(n - 2)}}
{{\text{2}}}}}{x^2} + \frac{1}
{6}{2^{ - \frac{{n(n - 1)}}
{{\text{2}}}}}{x^3} = 0,
\end{equation}
namely,
\begin{equation}{x^3} + 3 \cdot {2^{n - 1}}{x^2} + 3 \cdot {2^{2n - 2}}x + 3 \cdot {2^{3n - 5}} = 0.
\end{equation}
Obviously, the discriminant of the above cubic equation reads
\begin{equation}\Delta  = {2^{6n - 12}} > 0.
\end{equation}
Therefore, the cubic equation (4.148) has only one real root and a pair of imaginary roots, this contradicts to the assumption.

{\bf Step V.}$\;\;$  For any fixed integer $N>2$, we claim that the polynomial ${P_N}(x)$ has $N-2$ different real zeros denoted by
\begin{equation}{\{ {x_k}\} _{1 \le k \le N - 2}}\quad  {\rm with} \;\; {x_k} > {x_{k + 1}},
\end{equation}
and it holds that
\begin{equation}{\xi _{k + 1}} < {x_k} <  {\xi _{k - 1}} \quad \left(k = N - 2j - 1;\;\; j = 1,\cdots ,\left[ {\frac{{N - 1}}{2}} \right]\right)
\end{equation}
and
\begin{equation}{\xi _{N - 2}} < {x_{N - 2}}.
\end{equation}

It is easy to check that the statement is true for the case that $N=3,4,5,6,7$. We next prove this statement by induction. To do so, we firstly assume that the statement holds for integer $N-2\;(N>7)$, then we prove it also holds for integer $N$.

By assumption, the polynomial ${P_{N - 2}}(x)$ has $N-4$ different real zeros denoted by \begin{equation}{\{ {w_k}\} _{1 \le k \le N - 4}}\quad  {\rm with} \;\; \ {w_k} > {w_{k + 1}},\end{equation} and it holds that
\begin{equation}{\xi _{k + 1}} < {w_k} <  {\xi _{k - 1}}\quad \left(k = N - 2j - 1;\;\; j = 1,\cdots,\left[{\frac{{N - 3}}{2}} \right]\right),
\end{equation}
and
\begin{equation}{\xi _{N - 4}} < {w_{N - 4}}.
\end{equation}

For integer $k = N - 2j-3\;(j \in {\mathbb{N}_ + })$ with $2k - 1 > N$, we have
\begin{equation}|{w_k}| \le |{w_{N - 5}}| < |{w_{N - 4}}| < (N - 3){2^{N - 5}}.
\end{equation}
Hence
\begin{equation}{P_N}({w_k}){( - 1)^{N - 1}} = \frac{{|w_k^{N - 1}{\rm{|}}}}{{(N - 1)!{2^{\frac{{(N - 1)(N - 2)}}{2}}}}} - \frac{{|w_k^N{\rm{|}}}}{{N!{2^{\frac{{N(N - 1)}}{2}}}}} = \frac{{|w_k^N{\rm{|}}}}{{N!{2^{\frac{{N(N - 1)}}{2}}}}}\left(N{2^{N - 1}} + {w_k}\right) > 0.
\end{equation}
Combining Step II, we find that ${P_{N}}(x)$ has zeros in intervals
\begin{equation}( {\xi _{k + 1}},{w_k})\end{equation} and
\begin{equation}({w_k}, {\xi _{k - 1}}),\end{equation}
where $k = N - 2j-3\; (j \in {\mathbb{N}_+ })$ with $2k - 1 > N$. Then it follows from Steps II and III that ${P_{N}}(x)$ has zeros in intervals
\begin{equation}( {\xi _{N - 2}}, - 5N{2^{N - 6}})\end{equation} and
\begin{equation}( - 5N{2^{N - 6}}, {\xi _{N - 4}}).\end{equation}
By Step I, ${P_{N}}(x)$ has zeros in intervals
\begin{equation}( {\xi _{k}}, {\xi _{k-1}}),\end{equation}
where $k \in \left\{1,\cdots,\left[ {\frac{{N + 1}}{2}} \right] \right\}$. Thus, we may claim that ${P_{N}}(x)$ has at least $N-2$ different real zeros.

In fact, we may divide into two cases to prove it:

\vskip 3mm
\noindent{\bf Case 1: $N \equiv 0(\bmod 4)$ or $N \equiv 1(\bmod 4)$}

In the present situation, ${P_{N}}(x)$ has zeros in the following intervals:
${x_{N - 2}} \in ({\xi _{N - 2}}, - 5N{2^{N - 6}})$, ${x_{N - 3}} \in (-5N{2^{N-6}}, {\xi _{N-4}})$, ${x_{N - 4}} \in ({\xi _{N - 4}}, {w_{N - 5}})$,
${x_{N - 5}} \in ({w _{N - 5}}, {\xi _{N - 6}})$, $\cdots$,
${x_{\left[ {\frac{{N + 1}}{2}} \right] + 2}} \in \left({\xi_{\left[ {\frac{{N + 1}}{2}} \right] + 2}}, {w_{\left[ {\frac{{N + 1}}{2}} \right] + 1}}\right)$,
${x_{\left[ {\frac{{N + 1}}{2}} \right] + 1}} \in \left({w_{\left[ {\frac{{N + 1}}{2}} \right]+1}}, {\xi_{\left[ {\frac{{N + 1}}{2}} \right]}}\right)$,
${x_{\left[ {\frac{{N + 1}}{2}} \right]}} \in \left({\xi_{\left[ {\frac{{N + 1}}{2}} \right]}}, {\xi_{\left[ {\frac{{N + 1}}{2}} \right]-1}}\right)$,
$\cdots$, ${x_1} \in ({\xi_{1}}, {\xi_{0}})$.

\vskip 3mm
\noindent{\bf Case 2: $N \equiv 2(\bmod 4)$ or $N \equiv 3(\bmod 4)$}

In this case, ${P_{N}}(x)$ has zeros in the following intervals:
${x_{N - 2}} \in ({\xi _{N - 2}}, - 5N{2^{N - 6}})$,
${x_{N - 3}} \in (-5N{2^{N-6}}, {\xi _{N-4}})$,
${x_{N - 4}} \in ({\xi _{N - 4}}, {w_{N - 5}})$,
${x_{N - 5}} \in ({w _{N - 5}}, {\xi _{N - 6}})$,
$\cdots$,
${x_{\left[ {\frac{{N + 1}}{2}} \right] + 2}} \in \left({w_{\left[ {\frac{{N + 1}}{2}} \right] + 2}}, {\xi_{\left[ {\frac{{N + 1}}{2}} \right] + 1}}\right)$,
${x_{\left[ {\frac{{N + 1}}{2}} \right] + 1}} \in \left({\xi_{\left[ {\frac{{N + 1}}{2}} \right]+1}}, {\xi_{\left[ {\frac{{N + 1}}{2}} \right]}}\right)$,
${x_{\left[ {\frac{{N + 1}}{2}} \right]}} \in \left({\xi_{\left[ {\frac{{N + 1}}{2}} \right]}}, {\xi_{\left[ {\frac{{N + 1}}{2}} \right]-1}}\right)$,
$\cdots$
${x_1} \in ({\xi_{1}}, {\xi_{0}})$.

\vskip 4mm
On the other hand, in Step IV we proved that ${P_{N}}(x)$ has at most $N-2$ real zeros, hence  ${P_{N}}(x)$ has exactly $N-2$ different real zeros denoted by ${\{ {x_k}\} _{1 \leqslant k \leqslant N - 2}}$ with ${x_k} > {x_{k + 1}}$.
Moreover, the intervals above show that
\begin{equation}{\xi_{k+1}} < {x_k} < {\xi_{k-1}} \quad \left(k = N - 2j - 1;\; j = 1,\cdots,\left[ {\frac{{N - 1}}{2}} \right]\right)
\end{equation}
and
\begin{equation}{\xi_{N-2}} < {x_{N - 2}}.
\end{equation}
Thus, the statement has been established for any integer $N > 2$.

{\bf Step VI.}$\;\;$   Combining Steps I and V gives
\begin{equation}|{x_k}| < (k + 2){2^k}\quad (k = 1,2,\cdots,N - 2).
\end{equation}
On the other hand, the relationship between between zeros and coefficients of a polynomial yields
\begin{equation}|{z_N}{|^2} = \frac{{N!{2^{\frac{{N(N - 1)}}
{2}}}}}
{{\left| {\prod\limits_{n = 1}^{N - 2} {{x_n}} } \right|}} > \frac{{N!{2^{\frac{{N(N - 1)}}
{2}}}}}
{{\prod\limits_{n = 1}^{N - 2} {|{x_n}|} }} > \frac{{N!{2^{\frac{{N(N - 1)}}
{2}}}}}
{{\frac{1}
{2}N!{2^{\frac{{(N - 1)(N - 2)}}
{2}}}}} = {2^N}.
\end{equation}
Thus, we finish the proof of Proposition 4.3. $\quad\quad\quad\blacksquare$

\vskip 4mm

By Proposition 4.3, we can prove
\begin{Proposition}The solution $y=f(z)$ has no imaginary zeros.
\end{Proposition}
\noindent{\bf Proof.} For any given $R>0$, on the disc $|z| \leqslant R$ it holds that
\begin{equation}\left|\sum\limits_{n = 0}^\infty  {\frac{{{z^n}}}
{{n!{2^{\frac{{n(n - 1)}}
{2}}}}}} \right| \leqslant \sum\limits_{n = 0}^\infty  {\frac{{|z{|^n}}}
{{n!{2^{\frac{{n(n - 1)}}
{2}}}}}}  \leqslant \sum\limits_{n = 0}^\infty  {\frac{{{R^n}}}
{{n!{2^{\frac{{n(n - 1)}}
{2}}}}}}.
\end{equation}
This shows that the power series corresponding to the function $f(z)$ is uniformly convergent in any close disc. Obviously, the function $f(z)$ is analytic in $\mathbb{C}$. On the other hand, for the real number $R>0$ mentioned above, by Proposition 4.3, on the area
$$\mathbb{D} = \{ z:|z| \leqslant R\} \backslash \{ z:\:\;\operatorname{Im} z = 0, - R \leqslant \operatorname{Re} z \leqslant 0\;\}$$
${P_{N}}(x)$ has no zeros and does not identically equal to zero obviously for suitably large $N$.

Suppose that $f(z)$ has an isolated zero $\alpha$ in $\mathbb{D}$, then there exists a small disc denoted by
$$\Omega = \{ z:|z - \alpha | \leqslant r\}  \subset \mathbb{D}$$
such that $f(z)$  is analytic and has no zeros on $\partial \Omega$.
Therefore, by the argument principle we obtain
\begin{equation}\frac{1}
{{2\pi i}}\int\limits_{\partial \Omega} {\frac{{f'(z)}}
{{f(z)}}dz = 1}.
\end{equation}
Nevertheless, for large $N$ it holds that
\begin{equation}\frac{1}
{{2\pi i}}\int\limits_{\partial \Omega} {\frac{{{P_N}'(z)}}
{{{P_N}(z)}}dz = 0}.
\end{equation}
By the uniform convergence in any close disc, we have
\begin{equation}\mathop {\lim }\limits_{N \to \infty } \frac{1}
{{2\pi i}}\int\limits_{\partial \Omega} {\frac{{{P_N}'(z)}}
{{{P_N}(z)}}dz = } \frac{1}
{{2\pi i}}\int\limits_{\partial \Omega} {\frac{{f'(z)}}
{{f(z)}}dz}.
\end{equation}
This is a contradiction because of (4.168) and (4.169). Thus we prove that the function $f(z)$ has no zeros in the area $\mathbb{D}$.
Since $R>0$ is randomly chosen, we finish the proof of Proposition 4.4.
 $\quad\quad\quad\blacksquare$

\vskip 4mm

Theorem 3.1 comes from Propositions 4.2 and 4.4 immediately.

\begin{Remark} By Lemma 4.3, we can obtain a more precise estimate on the parameter ${\theta _n}$ in Theorem 3.1. In fact, it holds that
$$\frac{66}{25n} < {\theta _n} < \frac{2167}{789n}$$
for large $n$. Moreover, we have
\begin{equation}\lim_{n \to \infty}n{\theta _n} = C\approx 2.744.
\end{equation}
\end{Remark}

We next prove Theorem 3.2.

\vskip 3mm

\noindent{\bf Proof of Theorem 3.2.} Since ${\{ {x_n}\} _{n \geq 1}}$ is the set of zeros of $f(x)$ and $f(x)$ satisfies
\begin{equation}f'(x) = f\left( {\frac{x}{2}} \right),
\end{equation}
${\{ {2x_n}\} _{n \geq 1}}$ is the set of critical points  of the function $f(x)$. On the other hand, noting that every zero ${x_n}$ is simple and it holds that
\begin{equation}{x_{n + 1}} < 2{x_n} < {x_n},
\end{equation}
we have
\begin{equation}{( - 1)^n}f(2{x_n}) > 0,
\end{equation}
that is to say, the signs of extrema are alternate.

Direct calculation yields
\begin{equation}f(2{x_k}) = f\left( - (k + {\theta _k}){2^k}\right) = \sum\limits_{n = 0}^\infty {( - 1)^n} {\frac{{{{(k + {\theta _k})}^n}}}
{{n!}}} {2^{nk - n(n - 1)/2}}\mathop  = \limits^\Delta  \sum\limits_{n = 0}^\infty {( - 1)^n} {{u_n}}.
\end{equation}
Then
\begin{equation}\frac{{{u_{n + 1}}}}
{{{u_n}}} = \frac{{k + {\theta _k}}}
{{n + 1}}{2^{n - k}}.
\end{equation}
This implies that
\begin{equation}
{u_n} {\rm \;\; is\;\; increasing\;\; for\;\; }n<k;\quad {u_n}  {\rm \;\; is\;\; decreasing\;\; for\;\;}  n>k.
\end{equation}
Moreover, noting
\begin{equation}{u_{2k + 2}} = \frac{{{{(k + {\theta _k})}^{2k + 2}}}}
{{(2k + 2)!{2^{k + 1}}}} < \frac{{{{(k + {\theta _k})}^{2k + 2}}}}
{{{{\left( {\frac{{2k + 2}}
{e}} \right)}^{2k + 2}}{2^{k + 1}}}} < {\left( {\frac{e}
{{2\sqrt 2 }}} \right)^{2k + 2}} < 1,
\end{equation}
we have
\begin{equation}
|f(2{x_k})| < \sum\limits_{n = 0}^{2k + 1} {{u_n}}  + \left|\sum\limits_{n = 2k + 2}^\infty  {{u_n}{{( - 1)}^n}} \right|
< (2k + 2){u_k} + 1 = \frac{{2{{(k + 1)}^{k{\text{ + 1}}}}}}
{{k!}}{{\text{2}}^{\frac{{k(k + 1)}}
{2}}} + 1 < 2(k{\text{ + 1}}){e^{k + 1}}{{\text{2}}^{\frac{{k(k + 1)}}
{2}}}.
\end{equation}
Hence, for large $k$ it holds that
\begin{equation}{\log _2}|f(2{x_k})| < \frac{1}
{2}{k^2} + 2k.
\end{equation}
Again, noting
\begin{equation}|{\theta _{n + 1}} - {\theta _n}| < 1,
\end{equation}
we can choose a smooth function $\theta (t)$ defined on ${\mathbb{R}_ + }$ such that
\begin{equation}\theta (n) = {\theta _n},\quad 0 < \theta (t) < 1,\quad |\theta '(t)| < 1.
\end{equation}
Since the smooth function $g(t)=-(t + \theta (t)){2^t}$ is strictly monotonous for $t>0$,
we can make the following transformation
\begin{equation}x =  - (t + \theta (t)){2^t}\quad (0 < \theta (t) < 1).
\end{equation}
Then
\begin{equation}{\log _2}|x| = t + {\log _2}(t + \theta (t)).
\end{equation}
Hence for large $|x|$, it holds that
\begin{equation}{\log _2}|f(x)| < \frac{1}
{2}{t^2} + 2t < \frac{1}
{2}{\left( {{{\log }_2}|x|} \right)^2} + 2\log |x|.
\end{equation}
Thus, the proof of Theorem 3.2 is completed.  $\quad\quad\quad\blacksquare$

\section{Functional heat equation}

This section is devoted to the study on the following Cauchy problem for a nonlocal heat equation
\begin{equation}\left\{\aligned
&{u_t} - {u_{xx}}(t/2,x) = 0\quad (t \geq 0, x \in \mathbb{R}),\\
&u(0,x) = \varphi (x),\endaligned.\right.\end{equation}
where $\varphi (x)$ is the initial data.

\begin{Theorem} Suppose that the initial data $\varphi (x)$ satisfies
\begin{equation}\varphi (x) \in {L^1}(\mathbb{R}) \cap C(\mathbb{R}),
\end{equation} and its Fourier transform $\hat{\varphi} (\lambda )$ satisfies \begin{equation}\hat{\varphi} (\lambda ) \in {L^1}(\mathbb{R}).\end{equation}
Furthermore, for any fixed time $t \geq 0$, suppose that
\begin{equation}\hat{\varphi}(\lambda )f( - {\lambda ^2}t) \in {L^1}(\mathbb{R})\quad {\rm and}\quad {\lambda ^2}\hat{\varphi}(\lambda )f( - {\lambda ^2}t) \in {L^1}(\mathbb{R}).\end{equation}
Then there exists a real-value continuous solution $u=u(t,x)$ to the Cauchy problem (5.1). Moreover, assume $\hat{\varphi}(\lambda )$ has no imaginary zeros and $u(t,x)\in {L^1}(\mathbb{R})$ for any fixed $t \geq 0$, then the Fourier transform of the solution $u=u(t,x)$ has an infinite number of real zeros and no imaginary zeros.
\end{Theorem}

\noindent{\bf Proof.} By making the Fourier transform with respect to the space variable $x$, the Cauchy problem (5.1) becomes
\begin{equation}\left\{\aligned
&\frac{{d\tilde u(t,\lambda)}}
{{dt}} =  - {\lambda ^2}\tilde u(t/2,\lambda),\\
&\tilde u(0,\lambda) =\hat{\varphi} (\lambda)\endaligned \right.\end{equation}
Direct calculation leads to
\begin{equation}\tilde u(t,\lambda ) =\hat{\varphi} (\lambda )f( - {\lambda ^2}t).
\end{equation}

We now claim  the function
\begin{equation}u(t,x) = \frac{1}
{{2\pi }}\int_{ - \infty }^{ + \infty } {\hat{\varphi} (\lambda )f( - {\lambda ^2}t){e^{i\lambda x}}d\lambda }
\end{equation}
is a real-value continuous solution to the nonlocal heat equation.

In fact, since $\varphi (x) \in {L^1}(\mathbb{R}) \cap C(\mathbb{R})$ and its Fourier transform $\hat{\varphi} (\lambda ) \in {L^1}(\mathbb{R})$, we have
\begin{equation}u(0,x) = \frac{1}
{{2\pi }}\int_{ - \infty }^{ + \infty } {\hat{\varphi} (\lambda ){e^{i\lambda x}}d\lambda }  = \varphi (x).
\end{equation}
In the last equality we have made use of the inversion Fourier transformation.

Noting (5.4) and using the dominated convergence theorem yields
\begin{equation}{u_t}(t/2,x) = \frac{1}
{{2\pi }}\int_{ - \infty }^{ + \infty } { - {\lambda ^2}\hat{\varphi} (\lambda )f( - {\lambda ^2}t/2){e^{i\lambda x}}d\lambda }
\end{equation}
and
\begin{equation}{u_{xx}}(t/2,x) = \frac{1}
{{2\pi }}\int_{ - \infty }^{ + \infty } { - {\lambda ^2}\hat{\varphi} (\lambda )f( - {\lambda ^2}t/2){e^{i\lambda x}}d\lambda } .
\end{equation}
Hence, $u=u(t,x)$ solves the equation (1.8).
Clearly, $u(t,x)$, ${u_t}(t,x)$ and ${u_{xx}(t,x)}$ are continuous.

On the other hand, we claim that $u(t,x)$ is a real-value function.

In fact, the conjugate of $u(t,x)$ is
\begin{equation}\bar u(t,x) = \frac{1}
{{2\pi }}\int_{ - \infty }^{ + \infty } {\int_{ - \infty }^{ + \infty } {\varphi (\mu )} } {e^{i\lambda \mu }}d\mu f( - {\lambda ^2}t){e^{ - i\lambda x}}d\lambda  = {( {\check{\varphi}(\lambda )f( - {\lambda ^2}t)})^ \wedge }
\end{equation}
where $\bullet^{\hat{}}$ means the Fourier transform of $\bullet$, while $\bullet^{\check{}}$ means the inverse Fourier transform of $\bullet$. Since
\begin{equation}u(t,x) = {(\hat \varphi (\lambda )f( - {\lambda ^2}t))^ \vee },
\end{equation}
it suffices to prove
\begin{equation}( {\check{\varphi}(\lambda )f( - {\lambda ^2}t)})^ \wedge =(\hat \varphi (\lambda )f( - {\lambda ^2}t))^ \vee.
\end{equation}

For any fixed $t \geq 0$, let
$${f_n}(\lambda ) = f( - {\lambda ^2}t){\chi _{[ - n,n]}}(\lambda )\;\;(n \in {\mathbb{N}_ + }).$$
Clearly, ${f_n} \in {L^1}(\mathbb{R})$. Noting that
\begin{equation}\hat{f_n}=2\pi\check{f_n}.
\end{equation}
and using the convolution formula and the inversion  Fourier formula, we have
\begin{equation}
( {\check{\varphi}(\lambda )f_n(\lambda)})^ \wedge =(\hat \varphi (\lambda )f_n(\lambda))^ \vee.
\end{equation}
Let $n \to \infty$ and we obtain the desired equality (5.13) by the dominated convergence theorem.
This shows that $u=u(t,x)$ is a real-value continuous solution to the nonlocal heat equation (1.8).

Noting that $u(t,x)\in {L^1}(\mathbb{R})$ for any fixed $t \geq 0$, we have
\begin{equation}\int_{ - \infty }^{ + \infty } {u(t,x){e^{ - i\lambda x}}dx}  = \hat \varphi (\lambda )f( - {\lambda ^2}t) .
\end{equation}
Since we have proven that $f(x)$ has an infinite number of real(negative) zeros and no imaginary zeros, combining this matter and the assumption on $\hat \varphi (\lambda )$, we have completed the proof. $\quad\quad\quad\blacksquare$

\section{Conclusions and remarks}

In this paper we introduce and investigate a new kind of functional (including ordinary and evolutionary partial) differential equations with the form (1.2). Typical examples of this new kind of functional equations reads, e.g., the equations (1.3)-(1.4) or (1.7)-(1.9), etc.
This kind of new functional differential equations is a new and powerful tool to study some problems arising from both mathematics and physics, more applications in particular to differential geometry and fundamental physics can be expected. In the present paper, we explore new philosophy enjoyed these equations, in particular, for some typical examples, we prove the global existence of smooth solutions, analyze some interesting properties enjoyed by these solutions, and illustrate the differences between this new class of equations and the traditional one. On the other hand, we may replay the variable $t/2$ in (1.2) as $t/\alpha$ ($\alpha>1$ is a positive constant), and can carry out a similar discussion and obtain similar results.

As the end of this paper, we state the following conjecture:

\begin{Conjecture} For any $\alpha >1$, let $y = {f_\alpha }(x)$ be the solution to the Cauchy problem for the equation $y'(x) = y\left( {\frac{x}
{\alpha }} \right)$ with the initial data $y(0) = 1$. Let $\mathbb{A}$ be the set of the zeros of $y = {f_\alpha }(x)$.  Then $\mathbb{A}$ is countably infinite and $\mathbb{A} \subset \mathbb{R}$. Moreover, $\mathbb{A}$ can be denoted by ${\{ {x_n}\} _{n \geq 1}}$with ${x_n} > {x_{n + 1}}$, and it holds that
\begin{equation}{x_n} =  - (n + {\theta _n}){\alpha ^{n - 1}},
\end{equation}
where $\theta_n\in (0, M_\alpha)$ is a parameter depending on $n$ and $\alpha$, in which ${M_\alpha }$ is a positive constant only depending on $\alpha$ but independent of $n$. Moreover, $\{\theta_n\}$ satisfies
\begin{equation}\mathop {\lim }\limits_{n \to \infty } n{\theta_n} = \frac{\alpha}{3}(\ln{h(\alpha)})^{\prime},
\end{equation}
where $h(\alpha)$ is defined by
\begin{equation} h(\alpha )= \sum\limits_{k = 1}^\infty  {(2k - 1){{( - 1)}^{k - 1}}{\alpha ^{\frac{{k - {k^2}}}
{2}}}}.
\end{equation}
\end{Conjecture}
\begin{Remark}
In fact, in a similar manner, we can easily prove the conjecture for the case $\alpha \geq 2$. For the case $\alpha\in (1,2)$, it is worthy to study in future.
\end{Remark}

\vskip 5mm

\noindent{\Large {\bf Acknowledgements.}}  This work was supported in part by the NNSF of China (Grant Nos.: 11271323, 91330105) and the Zhejiang Provincial Natural Science Foundation of China (Grant No.: LZ13A010002).

\end{document}